\documentclass[12pt]{amsart}
\usepackage{amsmath,amscd,amssymb,amsthm,array}

\usepackage{mathtools}

\usepackage{hyperref}

\usepackage{tikz-cd}

\newcommand {\del}{{\partial}}

\usepackage{DLdef1}

\let\ssec\subsection

\renewcommand {\ssbegin}[2][*]
 {\refstepcounter{subsection}%
\if#1*
\addcontentsline{toc}{subsection}{\thesubsection.\hskip 1pc #2}%
\else
\addcontentsline{toc}{subsection}{\thesubsection.\hskip 1pc #2. #1}%
\fi
 \def \secno {\gdef \secno {}{\ssecfont
\thesubsection.\hskip 2ex}%
 }%
 \begin{#2}}

\renewcommand {\sssbegin}[2][*]
 {\refstepcounter{subsubsection}%\label{sss#1}
\if#1*
\addcontentsline{toc}{subsubsection}{\thesubsubsection.\hskip 1pc #2}%
\else
\addcontentsline{toc}{subsubsection}{\thesubsubsection.\hskip 1pc #2. #1}
\fi
 \def \secno {\gdef \secno {}{\ssecfont \thesubsubsection.\hskip 2ex}%
 }%
 \begin{#2}}

\renewcommand {\parbegin}[2][*]
 {\refstepcounter{paragraph}%\label{ssss#1}
\if#1*
\addcontentsline{toc}{paragraph}{\theparagraph.\hskip 1pc #2}%
\else
\addcontentsline{toc}{paragraph}{\theparagraph.\hskip 1pc #2. #1}
\fi
 \def \secno {\gdef \secno {}{\ssecfont \theparagraph.\hskip 2ex}%
 }%
 \begin{#2}}

 \makeatletter
\@namedef{subjclassname@2020}{\textup{2020} Mathematics Subject Classification}
\makeatother

\title{Gordan-Rankin-Cohen operators on superstrings}

\author[V.~Bovdi, D.~Leites]{Victor Bovdi${}^{a,*}$, Dimitry Leites${}^{b}$}
\address{
 ${}^*$The corresponding author\\
 ${}^a$Department of Mathematics, UAEU, Al Ain, UAE \\
 \email{vbovdi@gmail.com}\\
${}^b$Department of Mathematics Stockholm University, Albanov\"agen 28, SE-114 19 Stockholm
 Sweden\\
\email{dimleites@gmail.com}\\
}

\begin{document}

\begin{abstract} We distinguish two classifications of bidifferential operators: between  (A) spaces of modular forms and (B) spaces of weighted densities.

(A) The invariant under the projective action of $\text{SL}(2;\mathbb{Z})$   binary differential operators between spaces of modular forms of integer or half-integer weight on the 1-dimensional manifold were found  by Gordan (called transvectants), rediscovered and classified by Rankin and Cohen (called brackets), and, in still another context, by Janson and Peetre. The invariant under the algebraic supergroup $\text{OSp}(1|2; \mathbb{Z})$ super modular forms of integer and half-integer weight on $(1|1)$-dimensional superstrings with contact structure were introduced, bidifferential operators between them classified and further studied by Gieres-Theisen, Cohen-Manin-Zagier, and Gargoubi-Ovsienko.

(B) For any complex weights, we classify the analogs of Gordan-Rankin-Cohen (briefly: GRC) binary differential operators between spaces of weighted densities invariant under $\mathfrak{pgl}(2)$. For any complex weights, we classify the analogs of GRC-operators between spaces of weighted densities invariant under the Lie superalgebra $\mathfrak{osp}(1|2)$. In the case of $(1|1)$-dimensional superstring without any additional structure, we also classify the analogs of GRC-operators between spaces of any weighted densities invariant under the Lie superalgebra $\mathfrak{pgl}(1|2)$.
\end{abstract}



\maketitle

\markboth{{\itshape Victor Bovdi\textup{,} Dimitry Leites}}
{{\itshape Gordan-Rankin-Cohen operators on superstrings}}

%--------------------------------------------------------------------
%--------------------------------------------------------------------
%--------------------------------------------------------------------
%\thispagestyle{empty}
%\tableofcontents
%\setcounter{tocdepth}{2}
%%%%%%%%%%%%%%%%%%%%%%%%%%%%%%%%%%%%%
\section{Introduction} \label{Intro}
%%%%%%%%%%%%%%%%%%%%%%%%%%%%%%%%%%%%%
By \textit{invariant
operators} acting in the spaces of tensor
fields on a~ fixed manifold
$M^n$ of dimension~ $n$ we mean operators whose
form is the same in any %(curvilinear)
coordinate system.
The importance of such operators became manifest after discovery of
the relativity theory; for a~review of the literature devoted to solution of Veblen's problem (classify ``natural" operators), see \cite{Gr}. Eventually it became clear (to J.~Bernstein, see \cite{BeLe}) that speaking about \textit{differential} operators it suffices to consider a~domain $U\subset M^n$ and pass from invariance with respect to the group of diffeomorphisms to invariance with respect to the Lie algebra $\fvect(n)$ of vector fields with polynomial (or formal power series) coefficients since the problem is, actually, local. For the same reason, in proofs it suffices to consider only polynomial functions; the answers, however, are applicable to smooth and analytic functions.

The advent of supersymmetry theory as a~potential language of the future, currently conjectural, SUSY GUTS (Supersymmetric Grand Unified Theories of all fundamental forces) made it natural to extend the search of invariant operators in the spaces of tensor fields on supermanifolds, see \cite{GLS}, where the differential operators invariant with respect to Lie superalgebra $\fvect(n|m)$ of vector fields with polynomial coefficients on the $(n|m)$-dimensional superspace, and some of several its infinite-dimensional simple Lie subsuperalgebras, are described. Physicists call supermanifolds of superdimension $(1|N)$ \textit{superstrings} and usually consider them endowed with a~contact distribution, as was done in \cite{Gi, GTh, CMZ, GO}, where the superizations of modular forms were introduced and studied for $N=1$, together with invariant operators.

It is also natural to classify the operators invariant with respect to (maximal and simple) subalgebras of $\fvect(n|m)$. 
%, e.g., the Lie subsuperalgebra preserving a~contact structure.
For $(n|m)=(1|0)$, there is only one such subalgebra, namely $\fsl(2)$, see \cite{LSh}, and all tensor bundles with irreducible fibers are line bundles; their sections are called \textit{weighted densities}. Recall that on any finite-dimensional supermanifold, the elements of the space $\cF_w$ of weighted densities of \textit{weight} $w\in\Cee$ are expressions of the form $f(\vvol)^w$, where $f\in\cF:=\cF_0$ is a~function and $\vvol$ is a~volume element, on which the Lie derivative $L_X$ along the vector field $X$ with divergence $\Div X$ acts by the formula
\be\label{vvol}
L_X(f(\vvol)^w)=(X(f)+wf\Div X)(\vvol)^w.
\ee
Recall that a~holomorphic function $\varphi$ on the upper half-plane is called a~\textit{modular form of weight}~ $k$ (``typically a~positive integer", see \cite{W}) if $\SL(2;\Zee)$ transforms $\varphi$ as follows
\be\label{RC}
\varphi(x)\mapsto \varphi\left(\nfrac{ax+b}{cx+d}\right)(cx+d)^k\qquad \text{for any}\quad  \begin{pmatrix}a&b\\
c&d\end{pmatrix}\in\SL(2;\Zee),
\ee
and satisfies a~certain growth condition which is irrelevant for us in this paper. For integer and half-integer weights, the $\SL(2;\Zee)$-invariant differential operators between spaces of modular forms are called
\textit{Bol} operators in the \textbf{unary} case; the \textbf{binary} operators are called
\textit{Gordan transvectors} or \textit{Rankin--Cohen} operators, introduced and rediscovered in different settings in \cite{Gor, R, Co, Z, JP}.

Clearly,
any function $f$ on the complex line satisfying condition \eqref{RC} is a~coefficient of the density of weight~ $-\nfrac{k}{2}$. \textbf{Two natural problems arise: classify differential operators  between}
\begin{equation}\label{2pr}
\text{\begin{minipage}[c]{14cm}~
\begin{itemize}
\item[(A)]    spaces  of modular forms invariant under (the  subgroups of) $\SL(2;\Zee)$ or its superizations;\\
\item[(B)]   spaces of weighted densities invariant under (the subalgebras of) $\fvect(1|N)$.\\
\end{itemize}
\end{minipage}}
\end{equation}
Here, we consider Problem (B) on $(1|N)$-dimensional superstrings for $N=0, 1$: classify the differential operators between spaces of weighted densities invariant with respect to the maximal simple Lie subsuperalgebra $\fg\subset \fvect(1|N)$ and call them, somewhat abusing the terminology, also Bol (resp., Gordan-Rankin-Cohen, briefly: GRC) operators in the unary (resp., binary) case.

The problems (A) and (B) are related but not identical.  A~classification of  operators between spaces of modular forms invariant under the action of $\SL(2;\Zee)$ or its congruence subgroups differs from a~classification of $\fpgl(2)$-invariant operators between spaces of weighted densities: e.g., compare the degree-1 classical Rankin-Cohen bracket
\be\label{Z1}
\left(f, g\right)\longmapsto
\left(lf'g-kfg'\right)\in M_{k+l+2}\qquad \text{for any}\quad  f\in M_{k} \quad \text{and}\quad  g\in M_{l},
\ee
where $M_{k}$ is the space of modular forms of weight $k$, see \cite[Formula (9)]{Z1}, with the GRC-operator given by formula~\eqref{10} in which only weights $\mu_1=\mu_2=0$ are possible for $n=1$.

In \cite{BLS2}, it was shown that the description of Bol operators in the sense of Problem (B) on supermanifolds without any structure is a~wild problem, except for superstrings of superdimension $(1|1)$ --- unlike the superstrings with a~contact structure of any superdimension $(1|N)$, in which case the answer is given in \cite{BLS1}; some of these super Bol operators are particularly interesting. For a~classification of Bol operators in the sense of Problem (A), see \cite{Gi}.

Here, we classify (a)~ the GRC-operators on the line for any complex weights of the weighted densities, thus solving problem~ (B) on the line;
(b)~ the GRC-operators invariant under the only maximal simple Lie subsuperalgebra $\fosp(1|2)\subset \fk(1|1)$ on the $(1|1)$-dimensional superstring~ $\cU$ with a~contact structure --- continuation of \cite{Gi, GTh, CMZ, GO}; (c)~ the GRC-operators invariant under the only maximal simple Lie subsuperalgebra $\fpgl(2|1)\subset \fvect(1|1)$ on the $(1|1)$-dimensional superdomain~ $\cU$ without any additional structure, the case never considered before.

\section{The classical GRC-operators} \label{SScl}

\ssec{Classification of GRC-operators between spaces of modular forms}
Let $f(dx)^{\mu_1}$ (resp., $g(dx)^{\mu_2}$) be weighted densities of weight $\mu_1$ (resp., $\mu_2$) on a~1-dimensional manifold (string); then, their $n$th \textit{Rankin-Cohen bracket} (a.k.a. \textit{Gordan transvectant}) $[f,g]_n$ is the coefficient of $(dx)^{\mu_1+\mu_2+n}$ given for any $\mu_1,\mu_2\in\Cee$ and $n\in\Nee:=\{1, 2, \dots\}$ by the formula
\be\label{GRC}
\begin{array}{ll}
 [f,g]_{n}&:=\sum _{0\leq r\leq n}(-1)^{r}{\binom {\mu_1+n-1}{n-r}}{\binom {\mu_2+n-1}{r}}f^{(r)}g^{(n-r)},
\end{array}
\ee
where $\textstyle\binom {\mu}{k}:=\nfrac{\mu\cdot(\mu-1)\cdot\ldots\cdot(\mu-k+1)}{k!}$ for any $\mu\in\Cee$ and $k\in\Zee_+:=\{0, 1, \dots\}$, and $f^{(r)}:=(\partial_x)^r(f)$. The brackets \eqref{GRC} for positive integer $\mu_1$ and $\mu_2$ are ``the only universal bilinear combination of derivatives of $f$ and $g$ which goes from $M_{\mu_1}\otimes M_{\mu_2}$ to $M_{\mu_1+\mu_2+n}$", see \cite[p.26]{CMZ}.

\ssec{Classification of GRC-operators between spaces of weighted densities}
Having classified all bilinear differential operators between spaces of weighted densities invariant under any changes of variable (hence, under $\fpgl(2)$), see \cite{Gr}, Grozman discovered the Grozman operator
\[
Gz: \left(f(dx)^{-2/3},\ g(dx)^{-2/3}\right)\longmapsto
\left(2fg'''+3f'g''-3f''g'- 2f'''g\right)(dx)^{5/3}.
\]
Note that $Gz$ is also invariant under the congruence subgroup $\Gamma(3)\subset\SL(2;\Zee)$. But even speaking about integer weights, the unique order-1 differential operator from Grozman's classification
\be\label{10}
\left(f, g\right)\longmapsto
\left(af'g+bfg'\right)(dx)\qquad \text{for any}\quad \nfrac{a}{b}\in\mathbb{CP}^1
\ee
is also invariant under any changes of variable, hence is also a~GRC-operator. Let us classify, at last, all $\fpgl(2)$-invariant GRC-operators; Grozman prepared a~lot for this, see \cite{Gr}.

Let $\fg:=\fvect(1)=\fder\, \Cee[x]=\{f\partial\mid f\in\Cee[x]\}$, where $\partial:=\nfrac{d}{dx}$, with the grading $\deg x=1$.
All irreducible finite-dimensional $\fg_0=\fgl(1)$-modules are
1-dimensional; let $V^*_1$ and $V^*_2$ be such
modules. Let $v\in V^*_1$ and $w\in V^*_2$ be nonzero vectors of
weight $\mu_1$ and $\mu_2$, respectively, i.e.,
$$
(x\partial)v=\mu_1v, \qquad (x\partial)w=\mu_2w.
$$
We turn $V^*_1$ and $V^*_2$ into modules over $\fg_{\geq 0}:=\oplus_{i\geq 0}\ \fg_i$ by setting $X(V^*_1)=X(V^*_2)=0$ for any $X\in\fg_{>0}$, where $\fg_{> 0}:=\oplus_{i> 0}\ \fg_i$. Let $\partial '$ and $\partial ''$ be copies of $\partial$. Set
\[
\begin{split}
I(V^*_1, V^*_2):=I^{\fg}_{\fg_{\geq 0}}(V^*_1)\otimes I^{\fg}_{\fg_{\geq 0}}(V^*_2):&=(U(\fg)\otimes_{U(\fg_{\geq 0})}V^*_1)\otimes (U(\fg)\otimes_{U(\fg_{\geq 0})}V^*_2)\\
&\simeq (\Cee[\partial']\otimes V^*_1)\bigotimes (\Cee[\partial'']\otimes V^*_2).
\end{split}
\]
Let us describe the so-called \textit{singular} vectors in
$I(V^*_1, V^*_2)$, the ones killed by $X_+:=x^2\partial$. Let us
find all weight solutions of the equation
$X_+f_n=0$, where
$$
\textstyle
f_n=\mathop{\sum}\limits_{0\leq i\leq n}\frac{1}{i!(n-i)!}c_i(\partial ')^iv\otimes(\partial '')^{n-i}w\in I(V^*_1, V^*_2).
$$
Since on the line $\vvol=dx$, the restricted (sending polynomials to polynomials, not power series) dualization $(\cF_{-\mu_1})^*\simeq I(V^*_1)\simeq \Cee[\partial]\otimes V^*_1$ gives an interpretation of singular vectors $f_n$, see Theorem~\ref{ThGRC}, in terms of GRC-operators.

Observe that
\[
\begin{split}
X_+(\partial^i\otimes v)&=2i\partial^{i-1}\otimes(x\partial)v-i(i-1)
\partial^{i-1}\otimes v\\
&=i(2\mu_1-i+1)\partial^{i-1}v.
\end{split}
\]

Hence, the condition (in which the summands involving $(\del')^{-1}$ and $(\del'')^{-1}$ appear only with zero coefficients, so these summands make sense)
\[
\begin{split}
0=X_+f_n
=&\textstyle \mathop{\sum}\limits_{0\leq i\leq n}\frac{c_i}{i!(n-i)!} (i(2\mu_1-i+1)(\partial ')^{i-1}v\otimes
(\partial '')^{n-i}w\\
&\qquad +(n-i)(2\mu_2-n+i+1)(\partial ')^iv\otimes(\partial '')^{n-i-1}w)\\
=&\textstyle \mathop{\sum}\limits_{0\leq i\leq n-1}\frac{1}{i!(n-i)!}(c_{i+1}(2\mu_1-i)(\partial ')^iv\otimes(\partial '')^{n-i-1}w\\
&\qquad  +c_i(2\mu_2-n+i+1)(\partial ')^iv\otimes(\partial '')^{n-i-1}w)\\
\end{split}
\]
implies $n$ equations for $n+1$ unknowns
\be\label{chain}
(2\mu_2-n+i+1)c_i+(2\mu_1-i)c_{i+1}=0,\qquad \text{where}\quad i=0,\dots, n-1.
\ee

Denote by $M$ the $n\times (n+1)$ matrix whose elements are the coefficients of system \eqref{chain}. Clearly, $M$ can have non-zero elements only on the main diagonal --- the one that begins in the upper left corner --- and on the next one above it.
Denote the main diagonal by $D=(\alpha_0,\ldots, \alpha_{n-1})$, where $\alpha_i=2\mu_2-n+i+1$, and the one above it by $S=(\beta_0,\ldots, \beta_{n-1})$, where  $\beta_i=2\mu_1-i$. Obviously, at most one $\alpha_i$ (resp., $\beta_i$) can vanish. Indeed, if $\alpha_i=\alpha_j=0$ for some $i\not=j$, then $i=j$, which is impossible. Moreover, according to \eqref{chain}, if $\alpha_i=0$, then $\beta_i=0$ for the same $i$; and vice versa.

\underline{Case 1: $D$ has no zero components}. From \eqref{chain} we have:
\be\label{Vik01}
c_{i}=(-1)^{n-i}\textstyle\prod_{k=i}^{n-1}\frac{\beta_k}{\alpha_k}c_{n}=\prod_{k=i}^{n-1}\frac{k-2\mu_1}{2\mu_2-n+k+1}c_{n}, \qquad \text{where}\quad i=0,\ldots, n-1
\ee
so $\dim(Null(M))=1$ and $f_n$ depends on a~ parameter $c_{n}$; hence, a~unique solution up to a~factor.

\underline{Case 2: $D$ has a~zero component}. This case splits into several subcases.

\underline{Case 2a}. Let $n=1$ and $\alpha_0=\beta_0=0$. Then, $\dim(Null(M))=2$ and each $f_1$ depends on two  parameters $(c_{0}, c_{1})$; hence, up to a~factor, $f_1$ depends on $\nfrac{c_{0}}{c_{1}}\in\Cee\Pee^1$.

\underline{Case 2b}. Let $n\geq 2$ and $\alpha_0=\beta_0=0$. Then $\mu_1=0$, $\mu_2=\textstyle\frac{n-1}{2}$ and from \eqref{chain} we have:
\be\label{Vik02}
c_{l}=\textstyle\prod_{k=l}^{n-1}c_{n}=c_{n}^{n-l}, \qquad \text{where} \quad l=1,\ldots, n-1
\ee
so $\dim(Null(M))=2$ and each $c_l$ depends on two  parameters $(c_{0}, c_{n})$; hence, up to a~factor, $f_n$ depends on $\nfrac{c_{0}}{c_{n}}\in\Cee\Pee^1$.

\underline{Case 2c}. Let $n\geq 2$ and $\alpha_{n-1}=\beta_{n-1}=0$. Then $\mu_1=\textstyle\frac{n-1}{2}$, $\mu_2=0$ and from \eqref{chain} we have:
\be\label{Vik03}
c_{l}=\textstyle\prod_{k=l}^{n-2}c_{n-1}=c_{n-1}^{n-1-l}, \qquad \text{where} \quad l=0,\ldots, n-2
\ee
so $\dim(Null(M))=2$ and each $c_l$ depends on two  parameters $(c_{n-1}, c_{n})$; hence, up to a~factor, $f_n$ depends on $\nfrac{c_{n-1}}{c_{n}}\in\Cee\Pee^1$.

\underline{Case 2d}. Let $n\geq 2$ and $\alpha_{j}=\beta_{j}=0$ for $j\in\{ 1, \dots, n-2\}$. Then $\mu_1=\textstyle\frac{j}{2}$, $\mu_2=\textstyle\frac{n-j-1}{2}$ and from \eqref{chain} we have:
\be\label{Vik04}
\begin{split}
c_{i}&=\textstyle\prod_{k=i}^{j-1}\frac{k-2\mu_1}{2\mu_2-n+k+1}c_{j}=\prod_{k=i}^{j-1}c_{j}=c_{j}^{j-i}, \qquad \text{where} \quad i=0,\ldots, j-1\\
c_{l}&=\textstyle\prod_{k=l}^{n-1}\frac{k-2\mu_1}{2\mu_2-n+k+1}c_{n}=\prod_{k=l}^{n-1}c_{n}=c_{n}^{n-l}, \qquad \text{where} \quad l=j+1,\ldots, n-1\\
\end{split}
\ee
so $\dim(Null(M))=2$ and each $c_i$ and $c_l$ depends on two  parameters $(c_{j}, c_{n})$; hence, up to a~factor, $f_n$ depends on $\nfrac{c_{j}}{c_{n}}\in\Cee\Pee^1$.

\sssbegin{Theorem}\label{ThGRC} The classification of GRC-operators in terms of singular vectors is given by formulas \eqref{Vik01}--\eqref{Vik04} for the negatives of $\mu_1$ and $\mu_2$; i.e.,
\[
[f,g]_n=\textstyle\mathop{\sum}\limits_{0\leq i\leq n}\frac{1}{i!(n-i)!}c_if^{(i)}g^{(n-i)}.
\]
\end{Theorem}

\section{The GRC-operators on the $N=1$ superstring with contact structure} \label{SS11c}

Consider the $(1|1)$-dimensional superstring $\cM$ with coordinates $t$ (even) and $\theta$ (odd). Let a~distribution on $\cM$ be given by the 1-form $\alpha:=dt+\theta d\theta$. Let $\cF:=\Cee[t,\theta]$ be the space of functions and $\cF_\mu:=\cF\alpha^{\mu/2}$ the space of $\mu$-densities. The Lie superalgebra $\fk(1|1)$ is spanned by contact vector fields that preserve the distribution; they are of the form
\[
K_{f}:=\textstyle(2-\theta \partial_\theta)(f)\partial_t-(-1)^{p(f)}(\pderf{f}{\theta}
 - \theta\pderf{f}{t}) \partial_\theta\qquad \text{for any}\quad f\in\cF.
\]
The contact bracket of generating functions is given by the formula
\be\label{kb}
\{f,\, g\}:=\textstyle(2-\theta \partial_\theta)(f)\pderf{g}{t}-\pderf{f}{t}(2-\theta \partial_\theta)(g)-
(-1)^{p(f)}\pderf{f}{\theta} \pderf{g}{\theta}\qquad \text{for any}\quad f, g\in\cF.
\ee

\ssec{Classification of $\OSp(1|2;\Zee)$-invariant GRC-operators (\cite{Gi,GTh,CMZ, GO})}
%[\cite{GO}\footnote{This theorem and even comments after it were reproduced almost literally 11 years later by S.~Chouaibi and A.~Zbidi in J.~of Geometry and Physics, ignoring \cite{GO}.}]
\label{Ex} In the cited papers, it is assumed that $\deg \theta=\nfrac12$ and $\deg t=1$. Let $D_\theta=\partial_{\theta}-\theta \partial_t$, and $f^{(a)}:=(\partial_t)^{a}(f)$. Let $a\in\Zee_+$, and either $k\in\Nee:=\{1, 2, \dots\}$ or $k\in\Nee-\nfrac12:=\left\{\nfrac12, \nfrac32, \nfrac52, \dots\right\}$, let $\mu_1, \mu_2\in\Cee$ and 
\[
\gamma_{a,k}:=(-1)^{a+1} \binom{\mu_1+\lfloor k \rfloor}{\lfloor k \rfloor-a}\binom{\mu_2+\lfloor k \rfloor}{a}.
\]
The only invariant under the algebraic supergroup $\OSp(1|2;\Zee)$  bilinear differential operators between superspaces  of super modular forms $SM_{\mu}$ of integer and half-integer weight $\mu$  are the following ones
described in \cite{CMZ} together a proof of their uniqueness:
\[
{\cal J}_k: SM_{ \mu_1} \otimes SM_{\mu_2} \tto SM_{\mu_1+\mu_2+2k},
\]
explicitly given as follows for any $\varphi\alpha^{\mu_1/2}\in SM_{\mu_1}$ and $\psi\alpha^{\mu_2/2}\in SM_{\mu_2}$, where $\varphi$ and $\psi$ are functions.

(i) If $k\in\Nee$, then ${\cal J}_k$ is even; it is given by the formula
\[ 
\begin{array}{l}
{\cal J}_k (\varphi,\psi)=\sum\limits_{0\leq a\leq k -1} (-1)^{p(\varphi)} \gamma_{a,k-1} D_\theta( \varphi^{(a)})D_\theta (\psi^{(k -1-a)})- \sum\limits_{0\leq a\leq k} \gamma_{a,k} \varphi^{(a)} \psi^{(k -a)}.
\end{array}
\]

(ii) If $k\in\Nee-\nfrac12$, then ${\cal J}_k$ is odd; it is given by the formula

\[
{\cal J}_k (\varphi,\psi)=\sum_{0\leq a\leq \lfloor k \rfloor} \gamma_{a,k} \left ((-1)^{p(\varphi)} (\mu_1+a)\varphi^{(a)}D_\theta (\psi^{(\lfloor k \rfloor-a)})-(\mu_2+\lfloor k \rfloor-a)D_\theta (\varphi^{(a)})\psi^{(\lfloor k \rfloor-a)}\right ).
\]

\ssec{Classification of $\fosp(1|2)$-invariant GRC-operators} Let $\fg:=\fk(1|1)$, with the grading $\deg \theta=1$, $\deg t=2$, hence $\deg K_f=\deg f-2$ and $\fg_0=\Cee K_t$. The embedded subalgebra $i:\fosp(1|2)\tto\fk(1|1)$ is generated by $\nabla_-:=K_{\theta}$ and $\nabla_+:=K_{t\theta}$. We immediately see analogs of the operators~ \eqref{10}, due to the fact that $\Omega^1(\cM)/\cF\alpha\simeq \cF_1$. Indeed, let
\[
[df]:=df\pmod{\cF\alpha}\equiv \left(\nfrac{\partial f}{\partial \theta}-\theta\nfrac{\partial f}{\partial t}\right) d\theta=D_\theta(f)d\theta,
\]
where $df:=dt\nfrac{\partial f}{\partial t}+d\theta\nfrac{\partial f}{\partial \theta}$. The operator
\be\label{11}
\begin{split}
\left(f, g\right)\longmapsto
&\left(a[df]g+bf[dg]\right)\\
=&\left(aD_\theta(f)g+bfD_\theta(g)\right)d\theta\qquad \text{for any}\quad \nfrac{a}{b}\in\mathbb{CP}^1,\quad  f,g\in\cF
\end{split}
\ee
is invariant with respect to the whole $\fk(1|1)$, as was noted in \cite{GO}, hence is a~GRC-operator.

All irreducible finite-dimensional $\fg_0=\fgl(1)$-modules are
1-dimensional; let $V^*_1$ and $V^*_2$ be such
modules; let them be even. Let $v\in V^*_1$ and $w\in V^*_2$ be nonzero vectors of
weight $\mu_1$ and $\mu_2$, respectively, i.e.,
\[
K_t(v)=\mu_1v, \qquad K_t(w)=\mu_2w.
\]

 Setting $X(V^*_1)=X(V^*_2)=0$ for any $X\in\fg_{>0}$, where $\fg_{> 0}:=\oplus_{i> 0}\ \fg_i$ we turn $V^*_1$ and $V^*_2$ into modules over $\fg_{\geq 0}:=\oplus_{i\geq 0}\ \fg_i$. Let $K_\theta'$ and $K_\theta''$ be copies of $K_\theta$. Set $I(V^*_1, V^*_2):=I^{\fg}_{\fg_{\geq 0}}(V^*_1)\otimes I^{\fg}_{\fg_{\geq 0}}(V^*_2)$:
\[
\begin{split}
I(V^*_1, V^*_2):&=(U(\fg)\otimes_{U(\fg_{\geq 0})}V^*_1)\otimes (U(\fg)\otimes_{U(\fg_{\geq 0})}V^*_2)\\
&\simeq (\Cee[K_1, K_\theta']\otimes V^*_1)\bigotimes (\Cee[K_1, K_\theta'']\otimes V^*_2)\\
&{\simeq} (\Cee[K_\theta']\otimes V^*_1)\bigotimes (\Cee[K_\theta'']\otimes V^*_2)\qquad \text{since} \quad (K_\theta)^2=\nfrac12[K_\theta, K_\theta]=\nfrac12K_1.
\end{split}
\]
Let us describe the \textit{singular} vectors
$f_k\in I(V^*_1, V^*_2)$, the ones killed by $\nabla_+=K_{t\theta}$. Let us
find all weight (with respect to $K_t$) solutions of the equation
$\nabla_+f_k=0$, where
\[
\begin{split}
f_{2n}=\textstyle\mathop{\sum}\limits_{0\leq i\leq n}\nfrac{1}{i!(n-i)!}c_i(K_1')^iv&\otimes(K_1'')^{n-i}w\\
&+\textstyle\mathop{\sum}\limits_{0\leq i\leq n-1}e_iK_\theta(K_1')^iv\otimes K_\theta(K_1'')^{n-i-1}w;\\
f_{2n+1}=\textstyle\mathop{\sum}\limits_{0\leq i\leq n}a_iK_\theta(K_1')^iv&\otimes(K_1'')^{n-i}w\\
&+\textstyle\mathop{\sum}\limits_{0\leq i\leq n}\nfrac{1}{i!(n-i)!}b_i(K_1')^iv\otimes K_\theta(K_1'')^{n-i}w.\\
\end{split}
\]
We have $\vvol\simeq d\theta$, so the restricted dualization $(\cF_{-\mu_1})^*\simeq I(V^*_1)\simeq \Cee[K_\theta]\otimes V^*_1$ implies an interpretation of singular vectors $f_n$ in terms of GRC-operators. Note that the $\fk(1|1)$-invariant differential operators between spaces of tensor fields, in particular weighted densities, are polynomials in $D_\theta$ with constant coefficients, see \cite{Shch} and \cite[Theorem 3.2]{BLS1}; mind that $(D_\theta)^2=\partial_t$.

Note that $\{t\theta, 1\}=-2\theta$, $\{\theta, \theta\}=1$ and $\{t\theta, \theta\}=t$, see formula~\eqref{kb}. Then,
\[
\begin{split}
\nabla_+((K_1')^iv\otimes (K_1'')^{n-i}w)=&-2iK_\theta (K_1')^{i-1}v\otimes (K_1'')^{n-i}w\\
&-2(n-i)(K_1')^iv\otimes K_\theta (K_1'')^{n-i-1}w;\\
\end{split}
\]
\[
\begin{split}
\nabla_+(K_\theta (K_1')^iv\otimes& K_\theta (K_1'')^{n-i-1}w)\\
%=&K_t (K_1')^iv\otimes K_\theta (K_1'')^{n-i-1}w+\nfrac12 2i(K_1')^{i}v\otimes K_\theta (K_1'')^{n-i-1}w\\
&-K_\theta (K_1')^iv\otimes K_t (K_1'')^{n-i-1}w+\nfrac12\cdot  2(n-i-1)K_\theta (K_1')^iv\otimes (K_1'')^{n-i-1}w\\
=&(\mu_1-i)(K_1')^iv\otimes K_\theta (K_1'')^{n-i-1}w\\
&+ (-\mu_2+3(n-i-1)) K_\theta (K_1')^iv\otimes (K_1'')^{n-i-1}w;\\
&\\
\nabla_+(K_\theta (K_1')^iv\otimes& (K_1'')^{n-i}w)\\
=&(\mu_1-i)(K_1')^iv\otimes (K_1'')^{n-i}w+ 2(n-i)K_\theta(K_1')^iv\otimes K_\theta (K_1'')^{n-i-1}w;\\
&\\
\nabla_+((K_1')^iv\otimes& K_\theta (K_1'')^{n-i}w)\\
%=&-2iK_\theta (K_1')^{i-1}v\otimes K_\theta (K_1'')^{n-i}w+\\
%&+ (K_1')^iv\otimes K_t (K_1'')^{n-i}w-(-2)(n-i)(K_1')^iv\otimes K_\theta K_\theta(K_1'')^{n-i-1}w\\
%&=-2iK_\theta (K_1')^{i-1}v\otimes K_\theta (K_1'')^{n-i}w+ (\mu_2-2(n-i)+(n-i))(K_1')^iv\otimes (K_1'')^{n-i}w \\
=&-2iK_\theta (K_1')^{i-1}v\otimes K_\theta (K_1'')^{n-i}w+ (\mu_2-(n-i))(K_1')^iv\otimes (K_1'')^{n-i}w.\\
\end{split}
\]
The summands involving $(K_1')^{-1}$ and $(K_1'')^{-1}$ appear only with zero coefficients, so these summands make sense. We consider the two cases.

\underline{Case $\nabla_+(f_{2n})=0$}. The condition
\[
\begin{split}
0&=\nabla_+f_{2n}\\
&=\textstyle \mathop{\sum}\limits_{0\leq i\leq n}\frac{c_i}{i!(n-i)!}\left(-2iK_\theta (K_1')^{i-1}v\otimes (K_1'')^{n-i}w-2(n-i)(K_1')^iv\otimes K_\theta (K_1'')^{n-i-1}w\right)+\\
&\textstyle \mathop{\sum}\limits_{0\leq i\leq n-1}\frac{e_i}{i!(n-1-i)!}((\mu_1-i)(K_1')^iv\otimes K_\theta (K_1'')^{n-i-1}w\\
&\qquad\qquad + (-\mu_2+3(n+i+1)) K_\theta (K_1')^iv\otimes (K_1'')^{n-i-1}w)
\end{split}
\]
is equivalent to the following system
\be\label{Sect1}
\begin{split}
\textstyle\frac{(\mu_1-i)}{i!(n-1-i)!}e_i&=\textstyle\frac{2(n-i)}{i!(n-i)!}c_{i};\\
\textstyle\frac{(3(n+i+1) -\mu_2)}{i!(n-1-i)!}e_i&=\textstyle\frac{2(i+1)}{(i+1)!(n-i-1)!}c_{i+1}, \qquad \text{where}\quad i=0,\dots, n-1.
 \end{split}
\ee
Hence,
\be\label{Sect2}
\big(3(n+i+1) -\mu_2\big)e_{i}+\big((i+1)-\mu_1\big)e_{i+1}=0, \qquad \text{where}\quad i=0,\dots, n-1.
\ee

Similarly as we solved system \eqref{chain}, we denote by $M$ the $n\times (n+1)$ matrix whose elements are the coefficients of system \eqref{Sect2}. Clearly, $M$ can have non-zero elements only on the main diagonal $D=(\alpha_0,\ldots, \alpha_{n-1})$, where $\alpha_i=3(n+i+1) -\mu_2$,
and on the next one above it $S=(\beta_0,\ldots, \beta_{n-1})$, where  $\beta_i=(i+1)-\mu_1$. Obviously, at most one $\alpha_i$ (resp., $\beta_i$) can vanish. Indeed, if $\alpha_i=\alpha_j=0$ for some $i\not=j$, then $i=j$, which is impossible. Moreover, according to~ \eqref{Sect2}, if $\alpha_i=0$, then $\beta_i=0$ for the same $i$; and vice versa.

\underline{Case 1: $D$ has no zero components}. Then, from \eqref{Sect2} we have
\be\label{Hek00}
\begin{split}
e_{i}&=(-1)^{n-i}\textstyle\prod_{k=i}^{n-1}\frac{\beta_k}{\alpha_k}e_{n}=\prod_{k=i}^{n-1}\frac{\mu_1-(k+1)}{3(n+k+1) -\mu_2}e_{n};\\
c_{i}&=\textstyle\frac{\mu_1-i}{2}\prod_{k=i}^{n-1}\frac{\mu_1-(k+1)}{3(n+k+1) -\mu_2}e_{n}, \qquad \text{where}\quad i=0,\ldots, n-1.
\end{split}
\ee
Hence, $\dim(Null(M))=1$ and each of $e_i$ and $c_i$ depends on $e_{n}$; hence, up to a~factor, $f_1$ depends on $e_{n}$.

\underline{Case 2: $D$ has a~zero component}. Let $\alpha_j=\beta_j=0$ for some $j\in\{ 0, \dots, n-1\}$. Consider the following subcases:

\underline{Case 2(a).} Let $n=1$. Then $j=0$, so ${\mu_1}=1$ and $\mu_2=6$. Thus $\dim(Null(M))=2$ and depends on two  parameters $(e_{0}, e_{1})$; hence, up to a~factor,
\be\label{Hek01}
n=1,\quad \alpha_0=\beta_0=0,\quad {\mu_1}=1,\quad  \mu_2=6,\quad  f_1\quad \text{depends on}\quad \nfrac{e_{0}}{e_{1}}\in\Cee\Pee^1.
\ee

\underline{Case 2(b).} Let $n\geq 2$ and $\alpha_0=\beta_0=0$. Then $\mu_1=1$, $\mu_2=3(n+1)$ and from \eqref{Sect2} we have:
\be\label{Hek02}
\begin{split}
e_{l}&=\textstyle\prod_{k=l}^{n-1}\frac{\mu_1-(k+1)}{3(n+k+1) -\mu_2}e_{n}=(-\frac{1}{3})^{n-l}\textstyle\prod_{k=l}^{n-1}e_{n}=(-\frac{1}{3})^{n-l}e_{n}^{n-l}; \\
c_{l}&=\textstyle\frac{\mu_1-l}{2}(-\frac{1}{3})^{n-l}e_{n}^{n-l}, \qquad \text{where} \quad l=1,\ldots, n-1
\end{split}
\ee
so $\dim(Null(M))=2$ and each of the $e_i$ and $c_l$ depends on two  parameters $(e_{0}, e_{n})$; hence, up to a~factor, $f_n$ depends on $\nfrac{e_{0}}{e_{n}}\in\Cee\Pee^1$.

\underline{Case 2(c).} Let $n\geq 2$ and $\alpha_{n-1}=\beta_{n-1}=0$. Then $\mu_1=n$, $\mu_2=6n$ and from \eqref{Sect2} we have:
\be\label{Hek03}
\begin{split}
e_{l}&=\textstyle(-\frac{1}{3})^{n-l-1}\textstyle\prod_{k=l}^{n-2}\frac{n-(k+1)}{3(-n+k+1)}e_{n-1}=(-\frac{1}{3})^{n-l-1}e_{n-1}^{n-l-1};\\
c_{l}&=\textstyle\frac{\mu_1-l}{2}(-\frac{1}{3})^{n-l-1}e_{n-1}^{n-l-1},
 \quad \text{for} \quad l=0,\ldots, n-2
\end{split}\ee
so $\dim(Null(M))=2$ and each of the $e_l$ and $c_l$ depends on two  parameters $(e_{n-1}, e_{n})$; hence, up to a~factor, $f_n$ depends on $\nfrac{e_{n-1}}{e_{n}}\in\Cee\Pee^1$.

\underline{Case 2(d).} Let $n\geq 2$ and $\alpha_{j}=\beta_{j}=0$ for $j\in\{ 1, \dots, n-2\}$. Then $\mu_1=j+1$, $\mu_2=3(n+j+1)$ and from \eqref{Sect2} we have:
\be\label{Hek04}
\begin{split}
e_{i}&=\textstyle(-\frac{1}{3})^{j-i}e_{j}^{j-i},\quad c_{i}=\textstyle\frac{j-i+1}{2}(-\frac{1}{3})^{j-i}e_{j}^{j-i}, \quad \text{where} \quad i=0,\ldots, j-1;\\
e_{l}&=\textstyle(-\frac{1}{3})^{n-l}e_{n}^{n-l},\quad c_{l}=\textstyle\frac{j-l+1}{2}(-\frac{1}{3})^{n-l}e_{n}^{n-l}, \quad \text{where} \quad l=j+1,\ldots, n-1,
\end{split}
\ee
so $\dim(Null(M))=2$ and each of the $e_i$ and $c_l$ depends on two  parameters $(e_{j}, e_{n})$; hence, up to a~factor, $f_n$ depends on $\nfrac{e_{j}}{e_{n}}\in\Cee\Pee^1$.

 \underline{Case $\nabla_+(f_{2n+1})=0$}. From
\[
\begin{split}
0&=\nabla_+f_{2n+1}\\
&=\textstyle \mathop{\sum}\limits_{0\leq i\leq n}{a_i}\left((\mu_1-i)(K_1')^iv\otimes (K_1'')^{n-i}w+ 2(n-i)K_\theta(K_1')^iv\otimes K_\theta (K_1'')^{n-i-1}w\right)+\\
&\textstyle \mathop{\sum}\limits_{0\leq i\leq n-1}{b_i}\left(-2iK_\theta (K_1')^{i-1}v\otimes K_\theta (K_1'')^{n-i}w+ (\mu_2-n+i)(K_1')^iv\otimes (K_1'')^{n-i}w\right )
\end{split}
\]
it follows that
\be\label{Sect4}
\begin{split}
(i-\mu_1)a_i&=(\mu_2-n+i)b_i, \qquad\qquad  \text{where}\quad i=0,\dots, n-1\\
a_i&=\textstyle\frac{i+1}{n-i}b_{i+1}, \quad \qquad\qquad \qquad \text{where}\quad i=0,\dots, n-2\\
a_{n-1}&=a_n(\mu_1-n)=0.\\
 \end{split}
\ee

Consider the following cases:

\underline{Case 1. Let $\mu_2-n+i\not=0$ for all $i\in\{0,\dots, n-1\}$.} Then, for $i=n-1$ from the 1st equation of \eqref{Sect4} we get that $(\mu_2-1)b_{n-1}=0$ and $a_0=\cdots=a_{n-1}=0$ by the 1st and 2nd equations of \eqref{Sect4}. It follows that  $b_0=\cdots=b_{n-1}=0$ by the 1st equation of \eqref{Sect4} and $a_n(\mu_1-n)=0$. Consequently,
\be\label{Case1}
\text{if $\mu_2\not\in\{1,\dots, n\}$ and $\mu_1=n$, then there is a~unique (up to a~factor) solution.} 
\ee

\underline{Case 2. Let $\mu_2-n+k=0$ for some $k\in\{0,\dots, n-1\}$ and $\mu_1\not\in \{0,\dots, k-1\}$.} Similarly, as in the previous case, we obtain that
\[
\begin{split}
a_{k}=a_{k+1}&=\cdots=a_{n-1}=(\mu_1-n)a_n=0;\\
b_{k+1}&=b_{k+2}=\cdots=b_{n-1}=0.
\end{split}
\]
Clearly, if $b_k\not=0$ then $a_s\not=0$ for all $s\in\{0,\ldots, k-1\}$ by the 2nd equation of \eqref{Sect4}, so we have the following solution: 
\be\label{Case2}
\begin{split}
&b_i=\textstyle\frac{(i-\mu_1)(i+1)}{(\mu_2-n+i)(n-i)}b_{i+1}=\prod_{s=i}^{k-1}\frac{(s-\mu_1)(s+1)}{(s-k)(n-s)}b_{k}; \\
&a_i=\textstyle\frac{i+1}{n-i}\prod_{s=i+1}^{k-1}\frac{(s-\mu_1)(s+1)}{(s-k)(n-s)}b_{k}, \qquad\qquad  \text{where}\quad i=0,\dots, k-1\\
&a_{k}=a_{k+1}=\cdots=a_{n-1}=(\mu_1-n)a_n=0;\\
&b_{k+1}=b_{k+2}=\cdots=b_{n-1}=0,
\end{split}
\ee
which depends on two parameters $(b_{k}, a_{n})$ if $\mu_1=n$ and on one parameter $b_{k}$ otherwise; hence, $f_{2n+1}$ either depends on $\nfrac{b_{k}}{a_{n}}\in\Cee\Pee^1$ or is defined uniquely (up to a~factor), respectively.

\underline{Case 3. Let $\mu_2-n+k=0$ for some $k\in\{0,\dots, n-1\}$ and $\mu_1=t\in \{0,\dots, k-1\}$.} Then,
\[
(\mu_2-n+t)b_t=(t-k)b_t=0,
\]
where  $\mu_2-n+t\not=0$ because $t<k$. Hence, $b_t=0$ and similarly, as in two previous cases, we obtain that
\be\label{Case3}
\begin{split}
a_{0}&=\cdots=a_{t-1}=0; \qquad \qquad b_{0}=\cdots=b_{t}=0;\\
a_{k}&=\cdots=a_{n-1}=a_n=0;\quad b_{k+1}=\cdots=b_{n-1}=0;\\
b_i&=\textstyle\frac{(i-\mu_1)(i+1)}{(\mu_2-n+i)(n-i)}b_{i+1}=\prod_{s=i}^{k-1}\frac{(s-\mu_1)(s+1)}{(s-k)(n-s)}b_{k}; \\
a_i&=\textstyle\frac{i+1}{n-i}\prod_{s=i+1}^{k-1}\frac{(s-\mu_1)(s+1)}{(s-k)(n-s)}b_{k}, \qquad\qquad  \text{where}\quad i=t+1,\dots, k-1\\
\end{split}
\ee
which depends on two parameters $(a_{t}, b_{k})$; so, up to a~ factor, $f_{2n+1}$ depends on $\nfrac{a_{t}}{b_{k}}\in\Cee\Pee^1$.

\sssbegin{Theorem}\label{ThGRC2} The classification of GRC-operators in terms of singular vectors $f_k$ is given for the negatives of $\mu_1$ and $\mu_2$ by the formulas \eqref{Hek01}--\eqref{Hek04} for $k$ even and \eqref{Case1}--\eqref{Case3} for $k$ odd:
\[
\begin{split}
[f,g]_{2n}&=\textstyle\mathop{\sum}\limits_{0\leq i\leq n}\nfrac{1}{i!(n-i)!}c_if^{(i)}g^{(n-i)}\\
&\qquad +\textstyle\mathop{\sum}\limits_{0\leq i\leq n-1}\nfrac{1}{i!(n-1-i)!}e_iD_\theta(f^{(i)}) D_\theta(g^{(n-i-1)});\\
[f,g]_{2n+1}&=\textstyle\mathop{\sum}\limits_{0\leq i\leq n}\left(a_iD_\theta(f^{(i)})g^{(n-i)}+b_if^{(i)} D_\theta(g^{(n-i)})\right).
\end{split}
\] \end{Theorem}

\section{The GRC-problem on the general $N=1$ superstring} \label{SS11}
Recall that a~\textit{superdomain} $\cU$ is a~pair (a~domain $U$, the superalgebra $\cF$ of functions on $U$ with values in the Grassmann algebra of a~ vector space $V$); the morphisms of superdomains are in one-to-one correspondence with (parity preserving) automorphisms of the the superalgebra $\cF$, see \cite{Lo} containing further details. The superdimension of $\cU$ is the pair $(\dim U|\dim V)$.

Let the coordinates on the $(1|1)$-dimensional superdomain $\cU$ be $x$ (even) and $\xi$ (odd). Consider the standard $\Zee$-grading in $\fg:=\fvect(1|1)$, namely, we set $\deg x=\deg \xi=1$. Hence, $\fg_0\simeq\fgl(1|1)$. For a~basis of $\fg_{0}=\fgl(1|1)$ realized as a~subalgebra of $\fvect(1|1)$, we take (here $\del:=\del_x$ and $\delta:=\del_\xi$)
\[
X_-:=\xi\del, \quad H_1:=x\del, \quad H_2:=\xi\delta,\quad  X_+:=x\delta.
\]
The eigenvector with respect to $H_1$ and $H_2$ is called a~\textit{weight} vector, the weight vector annihilated by $X_+$ (resp., $X_-$) is called \textit{highest} (resp., \textit{lowest}). Let $V$ be a~$\fgl(1|1)$-module; setting $\fg_iV=0$ for $i>0$ we turn $V$ into a~module over $\fg_{\geq0}:=\oplus_{i\geq0}\fg_i$. The space of tensor fields of type $V$ is defined to be
\[
T(V):=\Hom_{U(\fg_{\geq0})}(U(\fg), V)\simeq\Cee[[x, \xi]]\otimes V.
\]

Hereafter, we work in terms of the spaces dual to $T(V)$, i.e.,
\[
I(V^*):=U(\fg)\otimes_{U(\fg_{\geq0})} V^*\simeq\Cee[\del, \delta]\otimes V^*.
\]

In what follows, the  ``singular vector" is the one annihilated by ${\fg_{>0}:=\oplus_{i>0}\fg_i}$. A generalization of \textbf{Veblen's bilinear problem} reduces to the task of classifying highest weight singular vectors in $I(V_1)\otimes I(V_2)$, where $V_1$ and $V_2$ are irreducible $\fgl(1|1)$-modules was solved, see a~review \cite{GLS}. The following problem is more cumbrous, but when feasible it can be solved by the same method.

\textbf{Problem (B) on a~$1|1$-dimensional superdomain}: Describe analogs of Gordan-Rankin-Cohen operators, i.e., $\fpgl(2|1)$-invariant operators
\[
T(V_1)\otimes T(V_2)\tto T(W)
\]
 between spaces of very special tensor fields --- spaces of weighted densities where all fibers $V_i$ and $W$ are 1-dimensional. In GRC-problem we consider the condition of singularity of the vector with respect to the part of $\fg_{>0}$ belonging to ${\fpgl(2|1)\subset \fvect(1|1)}$ embedded as a~graded subalgebra so that $\fpgl(2|1)_i=\fvect(1|1)_i$ for $i=-1$ and $0$ in the standard grading ($\deg x=\deg \xi=1$).

A~basis of $\fvect(1|1)_1$ can be formed by $s_x:=xE$ and $s_\xi:=\xi E=x\xi\del$, where $E=x\del+\xi\delta$, (these two vectors belong to $\fpgl(2|1)\subset \fvect(1|1)$) and two divergence-free elements which we do not need in this note: they do not belong to $\fpgl(2|1)$.

We denote $I(V_1):=\Cee[\del',\delta']\otimes V_1$ and $I(V_2):=\Cee[\del'',\delta'']\otimes V_2$, where $\del'$ and $\delta'$ (resp., $\del''$ and $\delta''$) are copies of $\del$ and $\delta$, as Grozman suggested in \cite{Gr}, and use the fact that
\[
I(V_1)\otimes I(V_2)\simeq \Cee[\del',\delta',\del'',\delta'']\otimes V_1\otimes V_2.
\]
Then, the singular vector of $I(V_1)\otimes I(V_2)$ on level $-n$ is of the form
\be\label{Singfn}
\begin{split}
f_n=\sum_{0\leq k\leq n} &(\del')^k (\del'')^{n-k} a_k \otimes r_k+\sum_{0\leq k\leq n-1} (\del')^{k} (\del'')^{n-k-1}\delta' \ b_k\otimes s_k \\[2mm]
&+
\sum_{0\leq k\leq n-1} (\del')^k(\del'')^{n-k-1} \delta'' \ c_k\otimes t_k+\sum_{0\leq k\leq n-2} (\del')^{k} (\del'')^{n-k-2} \delta'  \delta'' \ d_k\otimes u_k,
\end{split}
\ee
where $ a_k, b_k, c_k, d_k\in V_1$ and $r_k, s_k, t_k, u_k\in V_2$.

It suffices to confine ourselves to $\fg_0$-\textit{highest} singular vectors, the ones annihilated by $X_+$. Our task is to solve the system
\be\label{sys}
\begin{cases}
X_+(f_n)=0; \\
s_\xi(f_n)=0.
\end{cases}
\ee
Moreover, it suffices to confine ourselves to vectors
\be\label{f_n}
f_n\in \Cee[\del',\delta',\del'',\delta'']\otimes \text{irr}(V_1\otimes V_2),
\ee
where $\text{irr}(V_1\otimes V_2)$ is an irreducible $\fgl(1|1)$-component (sub- or quotient module) of $V_1\otimes V_2$, where the $V_i$ are irreducible $\fgl(1|1)$-modules.

%%%%%%%%%%%%%%%%%%%%%%%%%%%%%%%%%%%%%
\section{On $\fvect(1|1)$-modules induced from $\fgl(1|1)$-irreducibles} \label{Sprel}
%%%%%%%%%%%%%%%%%%%%%%%%%%%%%%%%%%%%%

 Let $\fg=\oplus_{i\geq -1}\ \fg_i$ be a~shorthand for $\fvect(1|1):=\fder\Cee[x, \xi]$ considered with the standard grading $\deg x=\deg\xi=1$. For a~basis of $\fg_{-1}$, we take $\del:=\del_x$ and $\delta:=\del_\xi$. Clearly, $\fg_0=\fgl(1|1)$.

\ssec{A description of irreducible $\fgl(1|1)$-modules with highest (or lowest) weight vector}
 Let $V:=W^{a;b}$ be an irreducible module over the commutative Lie algebra $\fh=\Span (H_1, H_2)$, and $v\in V$ a~non-zero element such that $H_1v=av$ and $H_2v=bv$ for some $a,b\in\Cee$. Since $\fh$ is commutative, $\dim V=1$. Let $\fb:=\fh\oplus \Cee X_+$. For definiteness, let $v$ be \textit{even}. We make $V$ into a~$\fb$-module by setting $X_+V=0$ and set
\[
\begin{split}
I_{\fb}^{\fg_{0}}(V):&=U(\fg_{0})\otimes_{U(\fb)}V\simeq U(\Cee X_-\oplus \fb)\otimes_{U(\fb)}V\\
&\simeq U(\Cee X_-)\otimes U(\fb)\otimes_{U(\fb)}V\simeq\Cee [X_-]\otimes V
\\
&
{\simeq}(\Cee X_-\oplus\Cee 1)\otimes V\qquad\qquad\qquad\qquad\qquad\qquad (\text{since}\quad X_-^2=\textstyle\frac12[X_-, X_-]=0)\\
& \simeq(\Cee X_-\otimes V)\oplus V.
\end{split}
\]

Clearly, $I_{\fb}^{\fg_{0}}(V)$ is naturally graded by powers of $X_-$; let $\deg V:=0$. Let us find out when there exists a~highest weight vector of degree $-1$ annihilated by $X_+$; we call such a~ vector \textit{singular}:
\[
X_+X_-v= Ev= (a+b)v=0\Rightarrow a+b=0.
\]

It is convenient to change notation and set
\[
H:=H_1-H_2\quad \text{and}\quad E:=H_1+ H_2.
\]
 Accordingly, set
\be\label{ab}
V^{\lambda;\mu}:=W^{a;b} \qquad \text{for}\quad  \lambda=a+b, \quad \mu=a-b.
\ee

Thus, if $\lambda\neq 0$, the module $M^{\lambda;\mu}:=I_{\fb}^{\fg_0}(V^{\lambda;\mu})$ is irreducible, while if $\lambda=0$, the module $M^{\lambda;\mu}$ is indecomposable, with a~submodule spanned by $X_-v$ and the quotient module spanned by $v$. The dual of this indecomposable module is also indecomposable with the submodule spanned by~ $v^*$. These indecomposable modules can be glued together in the same way as the adjoint $\fgl(1|1)$-module, whose maximal submodule is $\fsl(1|1):=\Span(X_-, E, X_+)$, is glued of 1-dimensional irreducibles. Note that $\fpgl(1|1):=\fgl(1|1)/\Cee E$. Let the actions of $X_-$ (resp., $X_+$) on the source of the arrow directed to the left (resp., right) and downwards be depicted as follows:

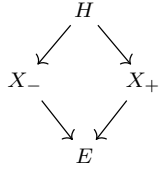
\[
\tiny
\begin{tikzcd}[column sep=tiny]
& H \ar[dl] \ar[dr] \\
X_- \ar[dr]&& X_+ \ar[dl] & \\
& E
&
&
\end{tikzcd}
\]
\normalsize

\subsection{The tensor product $\text{irr}(M^{\lambda;\, \mu})\otimes \text{irr}(M^{\sigma ; \, \rho })$} What are the irreducible $\fgl(1|1)$-modules ${M^{\lambda;\, \mu}\otimes M^{\sigma ; \, \rho}}$ is glued of? Recall that, by definition, the highest weight vectors $v\in M^{\lambda;\, \mu}$ and $w\in M^{\sigma ;\, \rho }$ are \textit{even}; if $v$ is odd, we denote the module it generates by $\Pi M^{\lambda;\, \mu}$. Let $\text{irr}(M^{\lambda;\,\mu})$ be the irreducible quotient of $M^{\lambda;\,\mu}$ modulo the maximal submodule; $\text{irr}(M^{\lambda;\,\mu})$ and $M^{\lambda;\, \mu}$ share the highest weight vector.

\sssec{The 4 cases}\label{4cases} Only the following cases are possible:

$(i)$ $\lambda=\sigma=0$. Clearly, $\text{irr}(M^{0;\, \mu})=V^{0;\, \mu}$ and $\text{irr}(M^{0; \, \rho })=V^{0;\, \rho }$; so $\dim V^{0;\, \mu}\otimes V^{0; \, \rho }=1$.

$(ii)$ $\lambda=0$ and $\sigma\neq 0$, then $V^{0;\,\mu}\otimes M^{\sigma; \, \rho }\simeq M^{\sigma;\, \rho +\mu}$ is irreducible. The case $\lambda\neq 0$ and $\sigma=0$ is analogous.

$(iii)$ $\lambda \sigma\neq 0$ and $\lambda +\sigma\neq 0$, then $M^{\lambda;\, \mu}\otimes M^{\sigma; \, \rho }$ is a~ direct sum of two irreducibles. Indeed, the vector
\be\label{subsp}
xX_- v\otimes w+yv\otimes X_- w\qquad \text{for some $x,y\in\Cee$}
\ee
is highest if
\[ 
xEv\otimes w+yv\otimes Ew=(x\lambda+y\sigma )v\otimes w=0.
\] 
The last condition is true only for
$x=-y\nfrac{\sigma }{\lambda}$.  Moreover,
\be\label{dir2}
M^{\lambda;\,\mu}\otimes M^{\sigma ; \, \rho }\simeq M^{\lambda+\sigma;\, \mu+\rho}\oplus \Pi M^{\lambda+\sigma;\, \mu+\rho-2},
\ee
where $\Pi$ is the reversal of parity functor.

$(iv)$ $\lambda \sigma\neq 0$ and $\lambda +\sigma=0$, then $M^{\lambda;\, \mu}\otimes M^{\sigma ; \, \rho }$ is the sum of 4 modules (two even and two odd), which are trivial as $\fsl(1|1)$-modules but non-trivial as $\fgl(1|1)$-modules, glued in the same way as the adjoint $\fgl(1|1)$-module is glued of irreducibles, see diagram \eqref{simi}.

Comparing with the structure of submodules in the adjoint $\fgl(1|1)$-module, we'd expect that if $\lambda +\sigma =0$, there are TWO highest weight vectors on ``level $-1$", not one.
Additionally, we'd expect that the vector $X_- v\otimes X_- w$ is a~highest weight vector, whereas (recall that $X_+$ and $X_-$ are odd, and the Sign Rule)
\[
\begin{split}
X_+(X_- v\otimes X_- w)&= Ev\otimes X_- w-X_- v\otimes Ew\\
&=
\lambda v\otimes X_- w-\sigma X_- v\otimes w= \lambda (v\otimes X_- w+X_- v\otimes w)\neq 0.
\end{split}
\]

However, there is no mistake, as I.~Shchepochkina explained.
Indeed, if $\lambda +\sigma =0$, there is a~1-dimensional invariant subspace
spanned by the vector $v\otimes X_- w+X_- v\otimes w$: both $X_-$ and $X_+$ annihilate it. The quotient modulo this subspace has 2 invariant 1-dimensional subspaces spanned by
$v\otimes w$ and $X_-v\otimes X_-w$, see formulas:
\[
\begin{array}{lll}
X_-:&v\otimes w&\mapsto X_-v\otimes w+v\otimes X_-w;\\
X_+:& X_-v\otimes w-v\otimes X_-w&\mapsto 2\lambda(v\otimes w);\\
X_-:&X_-v\otimes w-v\otimes X_-w&\mapsto -2X_-v\otimes X_-w;\\
X_+: &X_-v\otimes X_-w&\mapsto \lambda(v\otimes X_-w+X_-v\otimes w).\\
\end{array}
\]
The quotient modulo the direct sum of these subspaces is also
1-dimensional, so the structure of submodules is similar to that of the adjoint $\fgl(1|1)$-module (indicated are the basis vectors and their weights relative $H$ on the right of the vector, the label at the arrow is the coefficient, if distinct from $1$, the target is multiplied by):
\be\label{simi}
\begin{tikzcd}[column sep=tiny]
& {\scriptstyle X_-v\otimes w-v\otimes X_-w\ (\mu+\rho-2)} \ar[dl, "-2"] \ar[dr, "2\lambda"] \\
{\scriptstyle X_-v\otimes X_-w\ (\mu+\rho-4)}\ar[dr, "\lambda"]&&{\scriptstyle v\otimes w\ (\mu+\rho)} \ar[dl] & \\
& {\scriptstyle X_-v\otimes w+v\otimes X_-w\ (\mu+\rho-2)}
&
&
\end{tikzcd}
\ee

\ssec{Conclusion}\label{wild} The only feasible case to tackle in the search of singular vectors \eqref{Singfn}, \eqref{f_n} is the simplest Case $(i)$, see Subsection~\ref{4cases}, the classical setting of the GCR-problem. Perhaps, the method of generating functions, see \cite{KP}, can solve the other cases, but we do not know how to apply~ it.

In Case $(i)$, we consider the tensor product of the $\fvect(1|1)$-modules $I^{0;\, \mu}$ induced from the 1-dimensional $\fgl(1|1)$-modules of the form $V^{0;\, \mu}$. Recall that the module of weighted densities $\cF_a=\cF\vvol^a$ over the algebra of functions $\cF$ on the superstring $\cC^{1|1}$ is spanned by an element $\vvol^a$, where $a\in\Cee$, on which the Lie derivative $L_D$ along the vector field $D$ acts as multiplication by $a\Div(D)$, where the divergence is defined by the formula
\[
\textstyle
\Div(f\del+g\delta):=\frac{\del f}{\del x}+(-1)^{p(g)}\textstyle\frac{\del g}{\del {\xi}}\qquad\text{~~for any $f,g\in\cF$}.
\]

The restricted dualization $(\cF_{-\mu})^*\simeq I^{0;\,\mu}$ gives an interpretation of the results of Subsection~\ref{Sres}.
%%%%%%%%%%%%%%%%%%%%%%%%%%%%%%%%%%%%%

%%%%%%%%%%%%%%%%%%%%%%%%%%%%%%%%%%%%%

\ssec{Singular vectors in Case $(i)$} 
\label{Sres}
Let $H_1'$ (resp., $H_1''$) designate the operator $H_1$ acting on the 1st (resp., 2nd) factor in the tensor product $u:= v_1\otimes v_2\in V_1\otimes V_2$, where $V_1=V^{0;\mu_1}$ and $V_2= V^{0; \mu_2}$. Thus, we seek singular vectors of the form
\be\label{Singfn_2}
\begin{array}{ll}
f_n=\Big(\sum_{k=0}^{n}& A_k(\del')^k (\del'')^{n-k} +\sum_{k=0}^{n-1} B_k(\del')^{k} (\del'')^{n-k-1}\delta' \\[2mm]
&+
\sum_{k=0}^{n-1} C_k(\del')^k(\del'')^{n-k-1} \delta'' +\sum_{k=0}^{n-2} D_k(\del')^{k} (\del'')^{n-k-2} \delta'  \delta'' \Big)u.
\end{array}
\ee

To simplify the task, we use two observations due to A.~Lebedev. First, the summands of~ $f_n$ can be divided into the even (with coefficients $A$ and $D$) and odd (with coefficients $B$ and $C$), to be considered separately. Second, the even summands can be further subdivided: since all summands in $f_n$ are homogeneous with respect to the action of $H_2$, then their kernels are also direct sums of homogeneous subspaces. Hence, we can consider the summands with coefficients~ $A$ (of $H_2$-weight~ 0) and the summands with coefficients~  $D$ (of $H_2$-weight~ $-2$) separately.

Thus, the space of vectors of degree $-n$ can be subdivided into $3$ subspaces : 
\begin{itemize}

\item[(A)] with basis vectors $(\del')^k (\del'')^{n-k} u$, where $k=0,\dots, n$;

\item[(D)] with basis vectors $(\del')^k (\del'')^{n-k-2}\delta'\delta'' u$, where $k=0,\dots, n-2$;

\item[(BC)] with basis vectors $(\del')^k (\del'')^{n-k-1}\delta' u$, and $ (\del')^k (\del'')^{n-k-1}\delta'' u$, where $k=0,\dots, n-1$.
\end{itemize}
\underline{In case (A)}, we have
\[
X_+ (\del')^k (\del'')^{n-k} u = -k (\del')^{k-1} (\del'')^{n-k}\delta' u - (n-k)(\del')^k (\del'')^{n-k-1}\delta' u.
\]
The images of different basis elements are also different basis elements, and none of the images vanishes (except for $n=0$). Hence, the kernel of $X_+$ on the space of $A$-vectors is zero.

\underline{In case (D)}, we have
\be
\begin{split}
X_+ (\del')^k (\del'')^{n-k-2}\delta'\delta'' u =& 0;\\
s_\xi (\del')^k (\del'')^{n-k-2}\delta'\delta'' u =& (-k + \nfrac12\mu_1)(\del')^k (\del'')^{n-k-2}\delta'' u\\
& - (-(n-k-2) + \nfrac12\mu_2)(\del')^k (\del'')^{n-k-2}\delta' u.
\end{split}
\ee
The images of different basis elements are still different basis elements, but in this case, the image of the $k$-th basis element can vanish if $-k +\nfrac12 \mu_1$ and $-(n-k-2) + \nfrac12\mu_2$ vanish simultaneously.
Such $k$ exists if and only if $\mu_1$ and $\mu_2$ are non-negative even numbers and $\mu_1+\mu_2=2n-4$. We have
\begin{equation}\label{Ev}
\text{\begin{minipage}[c]{14cm} If $\mu_1=2k$ for some $k=0, 1, \dots, n-2$, then $\mu_2=2n-4-2k$ and $D_k\in\Cee$ is arbitrary whereas $D_i=0$ for $i\neq k$.\end{minipage}}
\end{equation}

\underline{In case (BC)}, we have
\[
\begin{split}
X_+ (\del')^k (\del'')^{n-k-1}\delta' u &= (n-k-1)(\del')^k (\del'')^{n-k-2}\delta'\delta'' u; \\
X_+ (\del')^k (\del'')^{n-k-1}\delta'' u &= -k(\del')^{k-1} (\del'')^{n-k-1}\delta'\delta'' u.
\end{split}
\]
This implies that the dimension of the kernel of $X_+$ on the space of BC-vectors is equal to $0|n+1$; and for a~basis of this space we can take
\be\label{els}
k(\del')^{k-1} (\del'')^{n-k}\delta' u + (n-k)(\del')^k (\del'')^{n-k-1}\delta'' u, \qquad \text{where $k=0,\dots, n$.}
\ee
Note that certain elements of this basis involve $(\del')^{-1}$ and $(\del'')^{-1}$ which appear only with zero coefficients, so these expressions make sense.

Now, let us see how $s_\xi$ acts on these elements \eqref{els}:
\[
\begin{split}
s_\xi (k(\del')^{k-1} (\del'')^{n-k}\delta' u &+ (n-k)(\del')^k (\del'')^{n-k-1}\delta'' u) \\
=& k(-(k-1)+\textstyle\frac12\mu_1)(\del')^{k-1} (\del'')^{n-k} u \\
&+ (n-k)(-(n-k-1)+\textstyle\frac12\mu_2)(\del')^k (\del'')^{n-k-1} u.
\end{split}
\]
In other words, if an element of the kernel of $X_+$ on the BC-space is of the form
\[
b=\sum_{0\leq k\leq n} b_k(k(\del')^{k-1} (\del'')^{n-k}\delta' u + (n-k)(\del')^k (\del'')^{n-k-1}\delta'' u),
\]
then the condition $s_\xi b = 0$ boils down to the conditions we obtain looking at the coefficient of $(\del')^k (\del'')^{n-k-1} u$ in $s_\xi b$):
\be\label{eqs}
\begin{split}
(n-k)(\mu_2-2&(n-k-1)) b_k \\&+ (k+1)(\mu_1-2k) b_{k+1} = 0, \qquad \text{where}\quad
k=0,\dots, n-1.
\end{split}
\ee

We solve this system in the same way as we solved system \eqref{chain}.
Denote by $M$ the $n\times (n+1)$ matrix whose elements are the coefficients of system \eqref{eqs}.
Clearly, $M$ can have non-zero elements only on the main diagonal and on the next one above it.  Let us denote the main diagonal by $D=(\alpha_0,\ldots, \alpha_{n-1})$, where $\alpha_i=(n-i)(\mu_2-2(n-i-1))$, and the one above it by $S=(\beta_0,\ldots, \beta_{n-1})$, where $\beta_i=(i+1)(\mu_1-2i)$. Obviously, at most one $\alpha_i$ (resp., $\beta_i$) can vanish. Indeed, if $\alpha_i=\alpha_j=0$ for some $i\not=j$, then $i=j$, which is impossible. Moreover, according to \eqref{eqs}, if $\alpha_i=0$, then $\beta_i=0$ for the same $i$; and vice versa.

We need to consider the following cases:

\underline{Case 1: $D$ has no zero component}. From \eqref{eqs} we have:
\be\label{Vikk3}
b_{i}=(-1)^{n-i}\textstyle\prod_{k=i}^{n-1}\frac{\beta_k}{\alpha_k}b_{n}=
\prod_{k=i}^{n-1}\frac{(k+1)(2k-\mu_1)}{(n-k)(\mu_2-2(n-k-1))}b_{n}, \qquad \text{where}\quad i=0,\ldots, n-1
\ee
so $\dim(Null(M))=1$ and $f_n$ depends on one parameter $b_{n}$; hence, there is a~unique solution of system \eqref{eqs}, up to a~factor.

\underline{Case 2: $D$ has a~zero component}. This case splits into several subcases.

\underline{Case 2a}. Let $n=1$ and $\alpha_0=\beta_0=0$. Then, $\mu_1=\mu_2=0$,
so $\dim(Null(M))=2$ and $f_1$ depends on two  parameters $(b_{0}, b_{1})$; hence, up to a~factor,
\be\label{2a}
n=1,\quad \alpha_0=\beta_0=\mu_1=\mu_2=0,\quad  f_1\quad \text{depends on}\quad \nfrac{b_{0}}{b_{1}}\in\Cee\Pee^1.
\ee

\underline{Case 2b}. Let $n\geq 2$ and $\alpha_0=\beta_0=0$. Then $\mu_1=0$, $\mu_2=2(n-1)$ and from \eqref{eqs} we have:
\be\label{Vikk4}
b_{l}=\textstyle\prod_{k=l}^{n-1}\textstyle\frac{(k+1)}{(n-k)}b_{n}, \qquad \text{where}\quad l=1,\ldots, n-1
\ee
so $\dim(Null(M))=2$ and each $b_l$ depends on two  parameters $(b_{0}, b_{n})$; hence, up to a~factor, $f_n$ depends on $\nfrac{b_{0}}{b_{n}}\in\Cee\Pee^1$.

\underline{Case 2c}. Let $n\geq 2$ and $\alpha_{n-1}=\beta_{n-1}=0$. Then $\mu_1=2(n-1)$, $\mu_2=0$ and from \eqref{eqs} we have:
\be\label{Vikk5}
b_{l}=\textstyle\prod_{k=i}^{n-2}\frac{(k+1)}{(k-n)}b_{n-1}, \qquad \text{where}\quad l=0,\ldots, n-2
\ee
so $\dim(Null(M))=2$ and each $b_l$ depends on two  parameters $(b_{n-1}, b_{n})$; hence, up to a~factor, $f_n$ depends on $\nfrac{b_{n-1}}{b_{n}}\in\Cee\Pee^1$.

\underline{Case 2d}. Let $n\geq 2$ and $\alpha_{j}=\beta_{j}=0$ for $j\in\{ 1, \dots, n-2\}$. Then $\mu_1=2j$, $\mu_2=2(n-j-1)$ and from \eqref{eqs} we have:
\be\label{Vikk6}
\begin{split}
b_{i}&=\textstyle\prod_{k=i}^{j-1}\frac{(k+1)(2k-\mu_1)}{(n-k)(\mu_2-2(n-k-1))}b_{j}
=\textstyle\prod_{k=i}^{j-1}\frac{(k+1)}{(n-k)}b_{j}, \qquad \text{where } \quad i=0,\ldots, j-1\\
b_{l}&=\textstyle\prod_{k=l}^{n-1}\frac{(k+1)(2k-\mu_1)}{(n-k)(\mu_2-2(n-k-1))}b_{n}
=\textstyle\prod_{k=l}^{n-1}\frac{(k+1)}{(n-k)}b_{n}, \qquad \text{where} \quad l=j+1,\ldots, n-1\\
\end{split}
\ee
so $\dim(Null(M))=2$ and each $b_i$ and $b_l$ depends on two  parameters $(b_{j}, b_{n})$; hence, up to a~factor, $f_n$ depends on $\nfrac{b_{j}}{b_{n}}\in\Cee\Pee^1$.

\sssbegin{Theorem}\label{ThMain} The dimension of the space of singular vectors $f_n$ is equal to
\begin{itemize}
\item[(i)] $1|1$ if $\mu_1$ and $\mu_2$ are non-negative even numbers and  $\mu_1+\mu_2=2n-4$;
\item[(ii)] $0|2$ if $\mu_1$ and $\mu_2$ are even integers between $0$ and $2n-2$ inclusive, and $\mu_1+\mu_2\geq 2n-2$;
\item[(iii)] $0|1$ otherwise.
\end{itemize}

So, up to a~non-zero constant factor, the odd singular vectors $f_n$ are uniquely defined in case $(iii)$ and case $(i)$ by the formula \eqref{Vikk3}. In case $(ii)$, there is a~$1$-parameter family of odd singular vectors given by formulas \eqref{2a}--\eqref{Vikk6}. In case $(i)$, there is also an even singular vector given by formula~\eqref{Ev}.
\end{Theorem}

\section{Discussion}

The classical GRC-operators (solution to Problem (A)) appear in number theory, representation theory, in functional analysis, and in theoretical physics; for useful references, see \cite{Z, Z1, Gi, GTh, KP, RTY, BSCK, BSCK1}. In this article, we solved Problem (B) using simple linear algebra, whereas solution of Problem (A) requires a complicated technique, see \cite{CMZ}. In the literature, speaking of Problem (A) Problem (B) was sometimes meant, especially in super setting.

\ssec{Open problems} In the case of 1-dimensional varieties over the ground field of characteristic $p>0$, all bilinear invariant operators between spaces of weighted densities are classified, see \cite{BL}, thus providing with examples of GRC-operators as well. To list all GRC-operators is a seemingly feasible open problem  if we restrict ourselves to tensor products of $\fsl(2)$-modules of dimension~ $\leq p$; for other open problems, see \cite{BL}.

In an unfinished draft of the description of GRC-operators between spaces of weighted densities on $(1|N)$-dimensional superstrings with a~ contact structure (by Bouarroudj et al.) many interesting new $\fk(1|N)$-invariant, hence GRC, operators are found. In particular, several incomplete for $N=2$  results of \cite{BLO,BBChM} dealing with Problem (B) are completed: for $N=2$, the weight is a pair of numbers and a~contact structure with odd time is possible.

In \cite{OR}, analogs of GRC-operators invariant under $\fo(n+2)\subset\fvect(n|0)$ between spaces of weighted densities are classified; it seems also feasible to classify GRC-operators invariant under other maximal finite-dimensional simple Lie sub(super)algebras of simple infinite-dimensional vectorial Lie (super)algebras (listed in \cite{LSh}).

However, if we generalize the setting of Problem (B) by considering binary operators between spaces of tensor fields with multidimensional fibers invariant with respect to $\fpgl(a+1|b)\subset \fvect(a|b)$, then Grozman's classification of the $\fvect(a)$-invariant operators, see \cite{Gr}, gives us  a~part of the answer. To tackle the problem in such generality is hardly feasible, see Grozman's proof and Subsection~\ref{wild}.

\subsection*{Acknowledgements} 
We are thankful to I.~Shchepochkina and A.~Lebedev for help.
%V.B. was supported by the grant ??

\def\eightit{\it}
\def\bib{\bf}
\bibliographystyle{amsalpha}

\end{document}